\documentclass{article}

\usepackage[T1]{fontenc}
\usepackage{textcomp}
\usepackage{ccfonts} %\usepackage{times}
\usepackage{euler}
\usepackage{amssymb}
\usepackage{amsmath}
\usepackage{amscd}
\usepackage{latexsym}
\usepackage{graphicx}   % to include and manipulate graphics
\usepackage[all]{xy}    % to draw diagrams
\usepackage{psfrag} % to add latex into eps files

\usepackage{epsfig}

\newtheorem{theorem}{Theorem}

\newtheorem{corollary}[theorem]{Corollary}

\newtheorem{lemma}[theorem]{Lemma}

\newtheorem{proposition}[theorem]{Proposition}
\newtheorem{assumption}[theorem]{Assumption}
\newtheorem{rema}{Remark}
\newenvironment{remark}{\begin{rema} \rm}{\end{rema}}
\newtheorem{exam}{Example}
\newenvironment{example}{\begin{exam} \rm}{\end{exam}}

\newenvironment{proof}[1][Proof]{\textbf{#1.} }{\ \rule{0.5em}{0.5em}}
\def\R{\mathbb{R}}

\def\E{\mathbb E}
\def\EXP{{\E}}

\def\PROB{{\mathbb P}}
\def\var{{\rm Var}}

\def\IND#1{{\mathbb I}_{{\left[ #1 \right]}}}

\def\prob{\PROB}
\def\wh{\widehat}
\def\ol{\overline}

\newcommand{\defeq}{\stackrel{\rm def}{=}}

\def\argmin{\mathop{\rm arg\, min}}

\def\sgn{\mathop{\rm sgn}}

\def\proof{\medskip \par \noindent{\sc proof.}\ }

\def\blackslug{\hbox{\hskip 1pt \vrule width 4pt height 8pt depth 1.5pt
\hskip 1pt}}
\def\qed{\quad\blackslug\lower 8.5pt\null\par\medskip}
\newcommand{\F}{{\cal F}}

\newcommand{\X}{{\cal X}}

\newcommand{\cR}{{\cal R}}
\def\N{{\cal N}}
\def\F{{\cal F}}

\def\Var{{\rm Var}}
\newcommand{\RR}{\mathbb{R}}
\newcommand{\NN}{\mathbb{N}}

\allowdisplaybreaks
\begin{document}
\title{Ranking and empirical minimization of $U$-statistics}
\author{St\'ephan Cl\'{e}men\c{c}on \\
MODALX - Universit\'{e} Paris X Nanterre\\
\&\\
Laboratoire de Probabilit\'{e}s et Mod\`{e}les Al\'{e}atoires\\
UMR CNRS 7599 - Universit\'{e}s Paris VI et Paris VII
\and G\'abor Lugosi
\thanks{The second author acknowledges support
by the Spanish Ministry of Science and Technology and
FEDER, grant BMF2003-03324
and by the PASCAL Network of
Excellence under EC grant no.\ 506778.
}
 \\
Departament d'Economia i Empresa\\
Universitat Pompeu Fabra
\and Nicolas Vayatis \\
Laboratoire de Probabilit\'{e}s et Mod\`{e}les Al\'{e}atoires\\
UMR CNRS 7599 - Universit\'{e}s Paris VI et Paris VII}

\maketitle

\bibliographystyle{plain}

\begin{abstract}
The problem of ranking/ordering instances, instead of simply
classifying them, has recently gained much attention in machine
learning. In this paper we formulate the \textit{ranking problem}
in a rigorous statistical framework. The goal is to learn a
ranking rule for deciding, among two instances, which one is
"better", with minimum ranking risk. Since the natural estimates
of the risk are of the form of a $U$-statistic, results of the
theory of $U$-processes are required for investigating the
consistency of empirical risk minimizers. We establish in
particular a tail inequality for degenerate $U$-processes, and
apply it for showing that fast rates of convergence may be
achieved under specific noise assumptions, just like in
classification. Convex risk minimization methods are also studied.
\end{abstract}

\section{Introduction}

Motivated by various applications including problems related to
document retrieval or credit-risk screening, the ranking problem
has received increasing attention both in the statistical and
machine learning literature. In the ranking problem one has to
compare two different observations and decide which one is
"better". For example, in document retrieval applications, one may
be concerned with comparing documents by degree of relevance for a
particular request, rather than simply classifying them as
relevant or not. Similarly, credit establishments collect and
manage large databases containing the socio-demographic and
credit-history characteristics of their clients to build a ranking
rule which aims at indicating reliability.

In this paper we define a statistical framework for studying such
ranking problems. The ranking problem defined here is closely
related to Stute's {\em conditional $U$-statistics}
\cite{Stu91,Stu94}. Indeed, Stute's results imply that certain
nonparametric estimates based on local $U$-statistics gives
universally consistent ranking rules. Our approach here is
different. Instead of local averages, we consider empirical
minimizers of $U$-statistics, more in the spirit of empirical risk
minimization popular in statistical learning theory, see, e.g.,
Vapnik and Chervonenkis \cite{VaCh74a}, Bartlett and Mendelson
\cite{BaMe06}, Bousquet, Boucheron, Lugosi \cite{BoBoLu05},
Koltchinskii \cite{Kol06}, Massart \cite{Mas06} for surveys and
recent development. The important feature of the ranking problem
is that natural estimates of the ranking risk involve
$U$-statistics. Therefore, the methodology is based on the theory
of $U$-processes, and the key tools involve maximal and
concentration inequalities, symmetrization tricks, and a
"contraction principle" for $U$-processes. For an excellent
account of the theory of $U$-statistics and $U$-processes we refer
to the monograph of de la Pe\~na and Gin\'e \cite{deGi99}.

Furthermore we provide a theoretical analysis of certain
nonparametric ranking  methods that are based on an empirical
minimization of convex cost functionals over convex sets of
scoring functions. The methods are inspired by boosting-, and
support vector machine-type algorithms for classification. The
main results of the paper prove universal consistency of properly
regularized versions of these methods, establish a novel tail
inequality for degenerate $U$-processes and, based on the latter
result, show that fast rates of convergence may be achieved for
empirical risk minimizers under suitable noise conditions.

We point out that under certain conditions, finding a good ranking
rule amounts to constructing a scoring function $s$. An important
special case is the bipartite ranking problem in which the
available instances in the data are labelled by binary labels
(good and bad). In this case the ranking criterion is closely
related to the so-called {\sc auc} (area under the "{\sc roc}"
curve) criterion (see the Appendix for more details).

The rest of the paper is organized as follows. In Section
\ref{preliminary}, the basic models and the two special cases of the
ranking problem we consider are introduced. Section \ref{empmin}
provides some basic uniform convergence and consistency results for
empirical risk minimizers. Section \ref{fast} contains the main
statistical results of the paper, establishing performance bounds
for empirical risk minimization for ranking problems. In Section
\ref{examples}, we describe the noise assumptions which guarantee
fast rates of convergence in particular cases. In Section
\ref{moment} a new exponential concentration inequality is
established for $U$-processes which serves as a main tool in our
analysis. In Section \ref{convex} we discuss convex risk
minimization for ranking problems, laying down a theoretical
framework for studying boosting and support vector machine-type
ranking methods. In the Appendix we summarize some basic
properties of $U$-statistics and highlight some connections of the
ranking problem defined here to properties of the so-called {\sc
roc} curve, appearing in related problems.

\section{The ranking problem}
\label{preliminary}

Let $(X,Y)$ be a pair of random variables taking values in $\X
\times \R$ where $\X$ is a measurable space. The random object $X$
models some observation and $Y$ its real-valued label. Let
$(X',Y')$ denote a pair of random variables identically
distributed with $(X,Y)$, and independent of it. Denote
\[
   Z= \frac{Y-Y'}{2}~.
\]
In the ranking problem one observes $X$ and $X'$ but not
 their labels $Y$ and $Y'$. We think about $X$ being
"better" than $X'$ if $Y>Y'$, that is, if $Z>0$. (The factor $1/2$
in the definition of $Z$ is not significant, it is merely here as
a convenient normalization.) The goal is to rank $X$ and $X'$ such
that the probability that the better ranked of them has a smaller
label is as small as possible. Formally, a {\em ranking rule} is a
function $r: \X\times \X \to \{-1,1\}$. If $r(x,x')=1$ then the
rule ranks $x$ higher than $x'$. The performance of a ranking rule
is measured by the {\it ranking risk}
\[
  L(r) = \PROB\{Z \cdot r(X,X') < 0\}~,
\]
that is, the probability that $r$ ranks two randomly drawn
instances incorrectly. Observe that in this formalization, the
ranking problem is equivalent to a binary classification problem
in which the sign of the random variable $Z$ is to be guessed
based upon the pair of observations $(X,X')$. Now it is easy to
determine the ranking rule with minimal risk. Introduce the
notation
\begin{align*}
\rho_+(X,X') & = \prob\{ Z >0 \mid X,X'\} \\
\rho_-(X,X') & = \prob\{ Z < 0 \mid X,X'\} \, .
\end{align*}
Then we have the following simple fact:

\begin{proposition}
Define
\[
r^*(x,x')=2\IND{\rho_+(x,x')\ge \rho_-(x,x')} -1
%  -\IND{\rho_+(x,x')<\rho_-(x,x')}
\]
and denote $L^*=L(r^*)=\EXP \{\min (\rho_+(X,X'),\rho_-(X,X'))\}$.
Then for any ranking rule $r$,
\[
   L^* \le L(r)~.
\]
\end{proposition}

\proof
Let $r$ be any ranking rule. Observe that, by conditioning first
on $(X, X')$, one may write
\[
L(r)=\EXP\left( \IND{r(X, X')=1} \rho_-(X,X') +\IND{r(X, X')=-1}
\rho_+(X,X')\right) .
\]
It is now easy to check that $L(r)$ is minimal for $r=r^*$.
\qed

Thus, $r^*$ minimizes the ranking risk over all possible ranking
rules. In the definition of $r^*$ ties are broken in favor of
$\rho_+$ but obviously if $\rho_+(x,x')= \rho_-(x,x')$, an
arbitrary value can be chosen for $r^*$ without altering its risk.

The purpose of this paper is to investigate the construction of
ranking rules of low risk based on training data. We assume that
$n$ independent, identically distributed copies of $(X,Y)$, are
available: $D_n=(X_1,Y_1),\ldots,(X_n,Y_n)$. Given a ranking rule
$r$, one may use the training data to estimate its risk
$L(r)=\PROB\{Z \cdot r(X,X') < 0\}$. The perhaps most natural
estimate is the {\em $U$-statistic}
\[
   L_n(r) = \frac{1}{n(n-1)}\sum_{i\neq j} \IND{Z_{i,j} \cdot r(X_i,X_j) <0}.
\]
In this paper we consider minimizers of the empirical estimate $L_n(r)$
over a class $\cR$ of ranking rules and study the performance of
such empirically selected ranking rules. Before discussing
empirical risk minimization for ranking, a few remarks are in order.

\begin{remark}
Note that the actual values of the $Y_i$'s are never used
in the ranking rules discussed in this paper. It is sufficient to
know the values of the $Z_{i,j}$, or, equivalently, the ordering
of the $Y_i$'s.
\end{remark}

\begin{remark}
{\sc (a more general framework.)} One may consider a
generalization of the setup described above. Instead of ranking
just two observations $X,X'$, one may be interested in ranking $m$
independent observations $X^{(1)},\ldots,X^{(m)}$. In this case
the value of a ranking function $r(X^{(1)},\ldots,X^{(m)})$ is a
permutation $ \pi$ of $\{1,\ldots,m\}$ and the goal is that $\pi$
should coincide with (or at least resemble to) the permutation
$\ol\pi$ for which $Y^{(\ol\pi(1))}\ge \cdots \ge
Y^{(\ol\pi(m))}$. Given a loss function $\ell$ that assigns a
number in $[0,1]$ to a pair of permutations, the ranking risk is
defined as
\[
   L(r) = \EXP \ell(r(X^{(1)},\ldots,X^{(m)}),\ol\pi)~.
\]
In this general case, natural estimates of $L(r)$ involve $m$-th
order $U$-statistics. Many of the results of this paper
may be extended, in a more or less
straightforward manner, to this general setup. In order to lighten
the notation and simplify the arguments, we restrict the
discussion to the case described
above, that is, to the case when $m=2$ and the loss function is
$\ell(\pi,\ol\pi) = \IND{\pi\neq \ol\pi}$.
\end{remark}

\begin{remark}
{\sc (ranking and scoring.)} In many interesting cases the ranking
problem may be reduced to finding an appropriate {\sl scoring
function}. These are the cases when the joint distribution of $X$
and $Y$ is such that there exists a function $s^*: \X \to \R$ such
that
\[
  r^*(x,x') = 1 \quad \mbox{if and only if} \quad s^*(x)\ge s^*(x')~.
\]
A function $s^*$ satisfying the assumption is called an {\sl
optimal scoring function}. Obviously, any strictly increasing
transformation of an optimal scoring function is also an optimal
scoring function. Below we describe some important special cases
when the ranking problem may be reduced to scoring.
\end{remark}

\begin{example}
{\sc (the bipartite ranking problem.)} In the bipartite ranking
problem the label $Y$ is binary, it takes values in $\{-1,1\}$.
Writing $\eta(x) = \PROB\{Y=1|X=x\}$, it is easy to see that the
Bayes ranking risk equals
\begin{align*}
L^* & =\EXP
\min \{\eta (X)(1-\eta(X')),\eta(X')(1-\eta(X))\} \\
&=\EXP \min \{\eta(X), \eta(X')\} -(\EXP \eta(X) )^2 \,
\end{align*}
and also,
\[
L^* =\Var\left(\frac{Y+1}{2}\right)-\frac{1}{2}\EXP\left| \eta
(X)-\eta (X')\right|~.
\]
In particular,
\[
L^*\leq \Var\left(\frac{Y+1}{2}\right)\leq 1/4
\]
where the equality $L^*=\Var\left(\frac{Y+1}{2}\right)$ holds when
$X$ and $Y$ are independent and the maximum is attained when $\eta
\equiv 1/2$. Observe that the difficulty of the bipartite ranking
problem depends on the concentration properties of the
distribution of $\eta (X)=\PROB (Y=1\mid X)$ through the
quantity $\EXP( \left|\eta (X)-\eta (X^{\prime })\right| ) $
which is a classical measure of concentration, known as  {\em
Gini's mean difference}.
For given $p=\EXP (\eta (X))$, Gini's mean difference ranges
from a minimum value of zero, when $\eta (X)\equiv p$, to a
maximum value of
$\frac{1}{2} p(1-p)$ in the case when
$\eta (X)= \left(Y+1\right)/2$. It is clear from the form of the
Bayes ranking rule that the optimal ranking rule is given by a
scoring function $s^*$ where $s^*$ is any strictly increasing
transformation of $\eta$. Then one may restrict the search to
 ranking rules defined by scoring functions $s$, that
is, ranking rules of form $r(x,x')=2\IND{s(x)\ge s(x')}-1$.
Writing $L(s)\defeq L(r)$, one has
\[
L(s)-L^*=\EXP \left( \left| \eta (X')-\eta(X)\right|
\IND{(s(X)-s(X'))(\eta (X)-\eta(X'))<0} \right) ~.
\]
We point out that the ranking risk in this case is closely related
to the {\sc auc} criterion which is a standard performance measure
in the bipartite setting (see \cite{FrIyScSi04} and Appendix 2). More
precisely, we have:
\begin{align*}
\mbox{{\sc auc}}(s) = & \PROB\left\{ s(X) \geq s(X^{\prime}) \mid Y=1, \,
Y^{\prime}=-1\right\} = 1-\frac{1}{2p(1-p)} L(s), \,
\end{align*}
where $p=\PROB\left( Y=1\right)$, so that maximizing the {\sc auc}
criterion boils down to minimizing the ranking error.
\end{example}

\begin{example}
{\sc (a regression model).} Assume now that $Y$ is real-valued and
the joint distribution of $X$ and $Y$ is such that $Y=m(X) +
\epsilon$ where $m(x)= \EXP(Y|X=x)$ is the regression function,
$\epsilon$ is independent of $X$ and
has a symmetric distribution around zero. Then clearly
the optimal ranking rule $r^*$ may be obtained by a scoring
function $s^*$ where $s^*$ may be taken as any strictly increasing
transformation of $m$.
\end{example}

\section{Empirical risk minimization}\label{empmin}

Based on the empirical estimate $L_n(r)$ of the risk $L(r)$ of a
ranking rule defined above, one may consider choosing a ranking
rule by minimizing the empirical risk over a class $\cR$ of
ranking rules $r:\X\times \X \to \{-1,1\}$. Define the empirical
risk minimizer, over $\cR$, by
\[
    r_n = \argmin_{r \in \cR} L_n(r)~.
\]
(Ties are broken in an arbitrary way.) In a "first-order"
approach, we may study the performance $L(r_n) = \PROB\{Z\cdot
r_n(X,X') < 0|D_n\}$ of the empirical risk minimizer by the
standard bound (see, e.g., \cite{DeGyLu95})
\begin{equation}
\label{supbound}
   L(r_n) - \inf_{r \in \cR} L(r) \le 2 \sup_{r\in \cR} |L_n(r)-L(r)|~.
\end{equation}
This inequality points out that bounding the performance of an
empirical minimizer of the ranking risk boils down to
investigating the properties of {\em $U$-processes}, that is,
suprema of $U$-statistics indexed by a class of ranking rules. For
a detailed and modern account of $U$-process theory we refer to
the book of de la Pe\~na and Gin\'e \cite{deGi99}. In a
first-order approach we basically reduce the problem
to the study of ordinary empirical
processes.

By using the simple Lemma \ref{utoemp} given in the Appendix,
we obtain the following:

\begin{proposition}
\label{firstorder}
Define the {\sl Rademacher average}
\[
   R_n = \sup_{r\in \cR} \frac{1}{\lfloor n/2 \rfloor}
       \left| \sum_{i=1}^{\lfloor n/2 \rfloor} \epsilon_i
    \IND{Z_{i,\lfloor n/2 \rfloor+i} r(X_i,X_{\lfloor n/2 \rfloor +i}) < 0}
       \right|
\]
where $\epsilon_1, ..., \epsilon_n$ are
i.i.d.\ Rademacher random variables (i.e., random symmetric sign variables).
Then for any convex nondecreasing function $\psi$,
\[
  \EXP \psi \left(L(r_n) - \inf_{r \in \cR} L(r) \right)
  \le \EXP \psi (4R_n)~.
\]
\end{proposition}

\proof
The inequality follows immediately from (\ref{supbound}),
Lemma \ref{utoemp} (see the Appendix),
and a standard symmetrization inequality,
see, e.g., Gin\'e and Zinn \cite{GiZi84}.
\qed

One may easily use this result to derive probabilistic performance
bounds for the empirical risk minimizer. For example, by taking
$\psi(x) = e^{\lambda x}$ for some $\lambda >0$, and using
 the bounded differences inequality (see
McDiarmid \cite{McD89}), we have
\begin{eqnarray*}
\lefteqn{
  \EXP \exp\left( \lambda( L(r_n) - \inf_{r \in \cR} L(r)) \right)  } \\*
& \le & \EXP \exp(4\lambda R_n)  \\*
& \le & \exp\left(4\lambda \EXP R_n
      + \frac{4\lambda^2}{(n-1)} \right)~.
\end{eqnarray*}
By using Markov's inequality and choosing $\lambda$ to minimize the
bound, we readily obtain:

\begin{corollary}
\label{firstorderbound}
Let $\delta >0$. With probability at least $1-\delta$,
\[
  L(r_n) - \inf_{r \in \cR} L(r) \le 4 \EXP R_n
    + 4 \sqrt{\frac{\ln (1/\delta)}{n-1}}~.
\]
\end{corollary}

The expected value of the Rademacher average $R_n$ may now be bounded by
standard methods, see,
e.g., Lugosi \cite{Lug02},
Boucheron, Bousquet, and Lugosi \cite{BoBoLu05}.
For example, if the class $\cR$ of indicator
functions has finite {\sc vc} dimension $V$, then
\[
\EXP R_n \le c\sqrt\frac{V}{n}
\]
for a universal constant $c$.

This result is similar to the one proved in the bipartite ranking
case by Agarwal, Graepel, Herbrich, Har-Peled, and Roth
\cite{AgGrHeHaRo05} with the restriction that their bound holds
conditionally on a label sequence. The analysis of
\cite{AgGrHeHaRo05} relies on a particular complexity measure
called rank-shatter coefficient but the core of the argument is
the same.

\medskip

The proposition above is convenient, simple, and, in a certain
sense, not improvable. However, it is well known from the theory
of statistical learning and empirical risk minimization for
classification that the bound (\ref{supbound}) is often quite
loose. In classification problems the looseness of such a
"first-order" approach is due to the fact that the variance of the
estimators of the risk is ignored and bounded uniformly by a
constant. Therefore, the main interest in considering
$U$-statistics precisely consists in the fact that they have
minimal variance among all unbiased estimators. However, the
reduced-variance property of $U$-statistics plays no role in the
above analysis of the ranking problem. Observe  that all upper
bounds obtained in this section remain true for an empirical risk
minimizer that, instead of using estimates based on
$U$-statistics, estimates the risk of a ranking rule by splitting
the data set into two halves and estimates $L(r)$ by
\[
\frac{1}{\lfloor n/2 \rfloor} \sum_{i=1}^{\lfloor n/2 \rfloor}
\IND{Z_{i,\lfloor n/2 \rfloor+i} \cdot r(X_i,X_{\lfloor n/2
\rfloor+i}) <0}~.
\]
Hence, in the previous study one loses the advantage of using
$U$-statistics. In Section \ref{fast}  it is shown that under
certain, not uncommon, circumstances significantly smaller risk
bounds are achievable. There it will have an essential importance
to use sharp exponential bounds for $U$-processes, involving
their reduced variance.

\section{Fast rates}\label{fast}

The main results of this paper show that the
bounds obtained in the previous section
 may be significantly improved under certain
conditions. It is well known (see, e.g., \S 5.2 in the survey
\cite{BoBoLu05} and the references therein) that tighter bounds
for the excess risk in the context of binary classification may be
obtained if one can control the variance of the excess risk by its
expected value. In classification this can be guaranteed under
certain "low-noise" conditions  (see Tsybakov \cite{Tsy02},
Massart and N\'ed\'elec \cite{MaNe06}, Koltchinskii \cite{Kol06}).

Next we examine possibilities of
obtaining such improved performance bounds for empirical ranking
risk minimization. The main message is that in the ranking problem
one also may obtain significantly improved bounds under some
conditions that are analogous to the low-noise conditions in the
classification problem, though quite different in nature.

Here we will greatly benefit from using $U$-statistics (as opposed
to splitting the sample) as the small variance of the
$U$-statistics used to estimate the ranking risk gives rise to
sharper bounds. The starting point of our analysis is
the Hoeffding decomposition of $U$-statistics (see Appendix 1).

Set first
\begin{equation*}
q_r((x,y),(x', y'))=\IND{(y-y') \cdot r(x,x')<0}-\IND{(y-y')\cdot
r^*(x,x')<0 }
\end{equation*}
and consider the following estimate of the {\em excess risk}
$\Lambda(r)=L(r)-L^*=\EXP q_r((X,Y),(X',Y'))$:
\begin{equation*}
\Lambda_n(r)=\frac{1}{n(n-1)}\sum_{i\neq j}
q_r((X_{i},Y_{i}),(X_{j},Y_{j})),
\end{equation*}
which is a $U$-statistic of degree $2$ with symmetric kernel
$q_r$. Clearly, the minimizer $r_n$ of the empirical ranking risk
$L_n(r)$ over $\cR$ also minimizes the empirical excess risk
$\Lambda_n(r)$. To study this minimizer, consider the Hoeffding
decomposition of $\Lambda_n(r)$:
\begin{equation*}
\Lambda_n(r) - \Lambda(r) = 2T_{n}(r)+W_{n}(r) \, ,
\end{equation*}
where
\[
T_{n}(r)=\frac{1}{n}\sum_{i=1}^{n} h_r(X_{i},Y_{i})
\]
is a sum of i.i.d. random variables with
\[
h_r(x,y) = \EXP q_r((x,y),(X',Y')) - \Lambda(r)
\]
and
\[
W_{n}(r)=\frac{1}{n(n-1)} \sum_{i\neq j}
\wh{h}_r((X_{i},Y_{i}),(X_{j},Y_{j}))
\]
is a {\sl degenerate} $U$-statistic with symmetric kernel
\[
\wh{h}_r((x,y),(x',y'))
   = q_r((x,y),(x',y'))-\Lambda(r)-h_r(x,y)-h_r(x',y') \, .
\]
In the analysis we show that the contribution of the  degenerate
part $W_n(r)$ of the $U$-statistic is negligible compared to that
of $T_n(r)$. This means that minimization of $\Lambda_n$ is
approximately equivalent to minimizing $T_n(r)$. But since
$T_n(r)$ is an average of i.i.d.\ random variables, this can be
studied by known techniques worked out for empirical risk
minimization.

The main tool for handling the degenerate part is a new
general moment inequality for $U$-processes that may be
interesting on its own right. This inequality is presented in
Section \ref{moment}. We mention here that for {\sc vc} classes
one may use an inequality of Arcones and Gin\'e \cite{ArGi94}.

It is well known from the theory of empirical risk minimization
(see Tsybakov \cite{Tsy02}, Bartlett and Mendelson \cite{BaMe06},
Koltchinskii \cite{Kol06}, Massart \cite{Mas06}), that, in order
to improve the rates of convergence (such as the bound
$O(\sqrt{V/n})$ obtained for {\sc vc} classes in Section
\ref{empmin}), it is necessary to impose some conditions on the
joint distribution of $(X,Y)$. In our case the key assumption
takes the following form:

\begin{assumption}
\label{varassumption}
There exist constants $c>0$ and $\alpha \in [0,1]$ such that
for all $r\in\cR$,
\[
\var(h_r(X,Y)) \leq c \,  \Lambda(r)^{\alpha}~.
\]
\end{assumption}

The improved rates of convergence will depend on the value of $\alpha$.
We will see in some examples that this assumption is satisfied
for a surprisingly large family of distributions, guaranteeing
improved rates of convergence.
For $\alpha=0$ the assumption is always satisfied and the corresponding
performance bound does not yield any improvement over those of
Section \ref{empmin}. However, we will see that in many natural
examples Assumption \ref{varassumption} is satisfied with
values of $\alpha$ close to one, providing significant improvements
in the rates of convergence.

Now we are prepared to state and prove the main result of the paper.
In order to state the result, we need to introduce some quantities
related to the class $\cR$. Let $\epsilon_1,\ldots,\epsilon_n$
be i.i.d.\ Rademacher random variables independent of
the $(X_i,Y_i)$. Let
\begin{align*}
Z_{\epsilon} & =  \sup_{r\in\cR} \left| \sum_{i,j} \epsilon_i
\epsilon_j \wh{h}_r((X_i,Y_i), (X_j,Y_j))  \right|~,\\
U_{\epsilon} &
= \sup_{r\in\cR} \sup_{\alpha : \|\alpha\|_2\le 1} \sum_{i,j} \epsilon_i  \alpha_j \wh{h}_r((X_i,Y_i),(X_j,Y_j)) ~,\\
M & =
\sup_{r\in\cR,k=1,\ldots,n} \left| \sum_{i=1}^n
\epsilon_i \wh{h}_r((X_i,Y_i),(X_k,Y_k)) \right|~.
\end{align*}
Introduce the "loss function"
\[
   \ell(r,(x,y)) = 2\EXP \IND{(y-Y) \cdot r(x,X)<0} - L(r)
\]
and define
\[
    \nu_n(r) = \frac{1}{n} \sum_{i=1}^n \ell(r,(X_i,Y_i))
       - L(r)~.
\]
(Observe that $\nu_n(r)$ has zero mean.)
Also, define the pseudo-distance
\[
d(r,r') =
  \left( \EXP \left( \EXP [\IND{r(X,X') \neq r'(X,X')} | X] \right)^2 \right)^{1/2}~.
\]
Let $\phi:[0,\infty) \to [0,\infty)$ be a nondecreasing function
such that $\phi(x)/x$ is nonincreasing and $\phi(1) \ge 1$ such
that for all $r\in \cR$,
\[
  \sqrt{n} \EXP \sup_{r' \in \R, d(r,r')\le \sigma}
 |\nu_n(r)-\nu_n(r')| \le \phi(\sigma)~.
\]

\begin{theorem}
\label{thm:fast}
Consider a minimizer $r_n$ of the  empirical ranking
risk $L_n(r)$ over a class $\cR$ of ranking rules and assume
Assumption \ref{varassumption}. Then there exists a
universal constant $C$ such that, with probability at least
$1-\delta$, the ranking risk of $r_n$ satisfies
\begin{eqnarray*}
  L(r_n) - L^* & \le & 2 \left(\inf_{r \in \cR} L(r) - L^* \right)  \\*
& & + C\left( \frac{\EXP Z_\epsilon}{n^2}
 + \frac{\EXP U_\epsilon  \sqrt{\log(1/\delta)}}{n^2}
 + \frac{\EXP M \log(1/\delta)}{n^2}
+    \frac{\log(1/\delta)}{n}  \right. \\*
& & + \left.
\rho^2 \log(1/\delta)
 \right)
\end{eqnarray*}
where  $\rho>0$ is the unique solution of the equation
\[
   \sqrt{n} \rho^2 = \phi(\rho^\alpha)~.
\]
\end{theorem}

The theorem provides a performance bound in terms of expected
values of certain Rademacher chaoses indexed by $\cR$ and local
properties of an ordinary empirical process. These quantities have
been thoroughly studied and well understood, and may be easily
bounded in many interesting cases. Below we will work out an
example when $\cR$ is a {\sc vc} class of indicator functions.

\proof
We consider the Hoeffding decomposition of the $U$-statistic
$\Lambda_n(r)$ that is minimized over $r\in \cR$.
The idea of the proof is to show that the
degenerate part $W_n(r)$ is of a smaller order and becomes
negligible compared to the part $T_n(r)$. Therefore,
$r_n$ is an approximate minimizer of $T_n(r)$ which can be
handled  by recent results on empirical risk minimization
when the empirical risk is defined as a simple sample average.

Let $A$ be the event on which
\[
\sup_{r \in \cR} |W_n(r)| \leq \kappa
\]
where
\[
\kappa =
C\left( \frac{\EXP Z_\epsilon}{n^2}
 + \frac{\EXP U_\epsilon  \sqrt{\log(1/\delta)}}{n^2}
 + \frac{\EXP M \log(1/\delta)}{n^2}
+    \frac{\log(1/\delta)}{n}  \right)
\]
for an appropriate constant $C$.
Then by Theorem \ref{thm:ineq}, $\PROB[A] \ge 1-\delta/2$.
By the Hoeffding
decomposition of the $U$-statistics $\Lambda_n(r)$
it is clear that, on $A$, $r_n$ is a $\rho$-minimizer
of
\[
\frac{2}{n} \sum_{i=1}^n \ell(r,(X_i,Y_i))
\]
over $r\in \cR$ in the sense that the value of this latter
quantity at its minimum is at most $\kappa$ smaller than at $r_n$.

Define $\tilde{r}_n$ as $r_n$ on $A$ and an arbitrary minimizer of
$(2/n) \sum_{i=1}^n \ell(r,(X_i,Y_i))$ on $A^c$. Then clearly, with
probability at least $1-\delta/2$, $L(r_n) = L(\tilde{r}_n)$ and
$\tilde{r}_n$ is a $\kappa$-minimizer of $(2/n) \sum_{i=1}^n
\ell(r,(X_i,Y_i))$. But then we may use Theorem 8.3 of
Massart \cite{Mas06} to bound the performance of $\tilde{r}_n$
which implies the theorem.
\qed

Observe that the only condition for the distribution is that the
variance of $h_r$ can be bounded in terms of $\Lambda(r)$. In
Section \ref{examples} we present examples in which Assumption
\ref{varassumption} is satisfied with $\alpha>0$. We will see
below that the value of $\alpha$ in this assumption determines the
magnitude of the last term which, in turn, dominates the
right-hand side (apart from the approximation error term).

The factor of $2$ in front of the approximation error term
$\inf_{r \in \cR} L(r) - L^*$ has no special meaning. It can be
replaced by any constant strictly greater than one at the price of
increasing the value of the constant $C$. Notice that in the bound
for $L(r_n)-L^*$ derived from Corollary \ref{firstorderbound}, the
approximation error appears with a factor of $1$. Thus, the
improvement of Theorem \ref{thm:fast} is only meaningful if
$\inf_{r \in \cR} L(r) - L^*$ does not dominate the other terms in
the bound. Ideally, the class $\cR$ should be chosen such that the
approximation error and the other terms in the bound are balanced.
If this was the case, the theorem would guarantee faster rates of
convergence. Based on the bounds presented here, one may design
penalized empirical minimizers of the ranking risk that select the
class $\cR$ from a collection of classes achieving this objective.
We do not give the details here, we just mention that the
techniques presented in Massart \cite{Mas06} and Koltchinskii
\cite{Kol06} may be used in a relatively straightforward manner to
derive such "oracle inequalities" for penalized empirical risk
minimization in the present framework.

In order to illustrate Theorem \ref{thm:fast}, we consider the
case when $\cR$ is a {\sc vc} class, that is, it has a finite
{\sc vc} dimension $V$.

\begin{corollary}
Consider the minimizer $r_n$ of the  empirical ranking
risk $L_n(r)$ over a class $\cR$ of ranking rules of finite {\sc vc}
dimension $V$ and assume
Assumption \ref{varassumption}. Then there exists a
universal constant $C$ such that, with probability at least
$1-\delta$, the ranking risk of $r_n$ satisfies
\[
  L(r_n) - L^*  \le  2 \left(\inf_{r \in \cR} L(r) - L^* \right)
 + C\left(
 \frac{V\log(n/\delta)}{n}
 \right)^{1/(2-\alpha)}
\]
\end{corollary}

\proof
In order to apply Theorem \ref{thm:fast}, we need  suitable
upper bounds for $\EXP Z_\epsilon$, $\EXP U_\epsilon$, $\EXP M$, and
$\rho$. To bound $\EXP Z_\epsilon$, observe that $Z_\epsilon$
is a Rademacher chaos indexed by $\cR$ for which Propositions 2.2 and 2.6
of Arcones and Gin\'e \cite{ArGi93} may be applied. In particular,
by using Haussler's \cite{Hau95} metric entropy bound for {\sc vc}
classes, it is easy to see that there exists a constant $C$ such that
\[
   \EXP Z_\epsilon \le CnV~.
\]
Similarly, $\EXP_\epsilon M$ is just an expected Rademacher average
that may be bounded by $C\sqrt{Vn}$ (see, e.g., \cite{BoBoLu05}).

Also, by the Cauchy-Schwarz inequality,
\begin{eqnarray*}
\EXP U_{\epsilon}^2 & \le & \EXP\sup_{r\in\cR} \sqrt{\sum_j \left(
\sum_i
\epsilon_i \wh{h}_r((X_i,Y_i),(X_j,Y_j)) \right)^2} \\*
& = & \EXP \sup_{r\in\cR} \left\{\sum_j \sum_i \wh{h}_r((X_i,Y_i),(X_j,Y_j))^2
\right.  \\*
& & \left. + \sum_j
\sum_{i,k} \epsilon_i \epsilon_k \wh{h}_r((X_i,Y_i),(X_j,Y_j))
    \wh{h}_r((X_j,Y_j),(X_k,Y_k))
 \right\}\\
 & \le & n^2 +
\EXP \sup_{r\in\cR} \sum_j
\sum_{i,k} \epsilon_i \epsilon_k \wh{h}_r((X_i,Y_i),(X_j,Y_j))
    \wh{h}_r((X_j,Y_j),(X_k,Y_k))~.
\end{eqnarray*}
Observe that the second term on the right-hand side is a Rademacher chaos
of order $2$ that can be handled similarly to $\EXP Z_\epsilon$.
By repeating the same argument, one obtains
\[
  \EXP U_{\epsilon}^2 \le n^2 + CVn^2
\]
Thus,
\[
\EXP(U_{\epsilon})\le
\sqrt{\EXP(U_{\epsilon}^2)}
\le Cn V^{1/2}~.
\]
This shows that the value of $\kappa$ defined in the proof of
Theorem \ref{thm:fast} is of the order of
$n^{-1}\left(V + \log(1/\delta)\right)$.
The main term in the bound of Theorem \ref{thm:fast} is $\rho^2$.
By mimicking the argument of Massart \cite[pp. 297--298]{Mas06},
we get
\[
C\left(
 \frac{V\log n}{n}
 \right)^{1/(2-\alpha)}
\]
as desired.
\qed

\section{Examples}
\label{examples}

\subsection{The bipartite ranking problem}

Next we derive a simple sufficient condition for achieving fast
rates of convergence for the bipartite ranking problem. Recall
that here it suffices to consider ranking rules of the form
$r(x,x')=2\IND{s(x)\ge s(x')}-1$ where $s$ is a scoring function.
With some abuse of notation we write $h_s$ for $h_r$.

\medskip

\noindent {\bf Noise assumption.} There exist constants $c>0$ and
$\alpha\in [0,1]$ such that for all $x\in\X$,
\begin{equation}
\label{NA1} \EXP_{X'}(\left| \eta (x)-\eta (X^{\prime })\right|
^{-\alpha})\leq c \, .
\end{equation}

\begin{proposition}
\label{noisecond} Under (\ref{NA1}), we have, for all $s\in\F$
\[
\var(h_s(X,Y)) \leq c \,  \Lambda(s)^\alpha \, .
\]
\end{proposition}

\proof
\begin{eqnarray*}
\lefteqn{ \var(h_s(X,Y)) } \\* & \leq & \EXP_{X}\left[ \left(
\EXP_{X'} (\IND{(s(X)-s(X^{\prime }))(\eta (X)-\eta
(X'))<0})\right) ^{2}\right]  \\* & \leq & \EXP_{X}\left[
\EXP_{X'}\left(\IND{(s(X)-s(X'))(\eta (X)-\eta (X'))<0}\left| \eta
(X)-\eta (X')\right|^{\alpha} \right) \right.  \\* & & \times
\left. \left( \EXP_{X'}(\left| \eta (X)-\eta (X')\right|
^{-\alpha})\right) \right] \\*
& & \quad \mbox{(by the Cauchy-Schwarz inequality)}  \\
& \le & c \left(\EXP_{X} \EXP_{X'} \left(\IND{(s(X)-s(X'))
(\eta(X)-\eta (X'))<0} \left| \eta (X)-\eta (X')\right| \right)
\right)^{\alpha}
\\*
& & \quad \mbox{(by Jensen's inequality and the noise assumption)}  \\
& = & c \Lambda(s)^{\alpha}~.
\end{eqnarray*}
\qed

%% \begin{remark}
%% Note that if we had considered the statistic based on a splitted sample, the same argument applied to $\var(q_s((X,Y),(X',Y'))$ would have given
%% \begin{multline*}
%% \var(q_s) \leq \EXP_{X}\left[ \EXP_{X'} \IND{s(X) - s(X')) \cdot (\eta (X)-\eta (X'))<0})\right]  \\
%% \leq \EXP_{X}\left[
%% \left( \EXP_{X'}(\IND{(s(X)-s(X'))(\eta
%% (X)-\eta (X'))<0})\left| \eta (X)-\eta (X')\right|
%% \right) ^{1/2}
%% \left( \EXP_{X'}(\left| \eta (X)-\eta (X')\right| ^{-1})\right) ^{1/2}
%% \right] \, .
%% \end{multline*}

%% Applying again Cauchy-Schwarz inequality, we  get $\var(h_s) \leq \sqrt{c \Lambda(s)}$ which eventually leads to a  slower rate (in $1/\sqrt{n}$).
%% \end{remark}

Condition (\ref{NA1}) is satisfied under quite general
circumstances. If $\alpha=0$ then clearly the condition poses no
restriction, but also no improvement is achieved in the rates of
convergence. On the other hand, at the other extreme, when
$\alpha=1$, the condition is quite restrictive as it excludes
$\eta$ to be differentiable, for example, if $X$ has a uniform
distribution over $[0,1]$. However, interestingly, for any
$\alpha<1$, it poses quite mild restrictions as it is highlighted
in the following example:

\begin{corollary}
Consider the bipartite ranking problem and assume that
$\eta(x)=\PROB\{Y=1|X=x\}$ is such that the random variable
$\eta(X)$ has an absolutely continuous distribution on $[0,1]$
with a density bounded by $B$. Then for any $\epsilon>0$,
\[
\forall x\in\X, \quad  \EXP_{X'}(\left| \eta (x)-\eta (X^{\prime
})\right|^{-1+\epsilon}) \leq \frac{2B}{\epsilon}
\]
and therefore, by Propositions \ref{fast} and \ref{noisecond},
there is a constant $C$ such that for
every $\delta,\epsilon \in (0,1)$,  the
excess ranking risk of the empirical minimizer $r_n$ satisfies,
with probability at least $1-\delta$,
\[
L(r_n) - L^* \leq
2\left( \inf_{r\in \cR} L(r) - L^* \right) +
 CB\epsilon^{-1} \left( \frac{V\log
(n/\delta)}{n} \right)^{1/(1+\epsilon)} \, .
\]
\end{corollary}

\proof
The corollary follows simply by checking that (\ref{NA1}) is
satisfied for any $\alpha=1-\epsilon<1$. Denoting the density of
$\eta(X)$ by $f$, we have
\begin{eqnarray*}
\EXP_{X'}(\left| \eta (x)-\eta (X')\right|^{-\alpha}) & = &
\int_0^1 \frac{1}{|\eta(x)-u|^\alpha} f(u)du  \\* & \le & B
\int_0^1 \frac{1}{|\eta(x)-u|^\alpha} du  \\* & = &
B\frac{\eta(x)^{1-\alpha} + (1-\eta(x))^{1-\alpha}}{1-\alpha} \le
\frac{2B}{1-\alpha}~.
\end{eqnarray*}
\qed

The condition (\ref{NA1}) of the corollary requires that the
distribution of $\eta(X)$ is sufficiently spread out, for example
it cannot have atoms or infinite peaks in its density. Under such
a condition a rate of convergence of the order of
$n^{-1+\epsilon}$ is achievable for any $\epsilon>0$.
\begin{remark}
Note that we crucially used the reduced variance of the
U-statistic $L(r_n)$ to derive fast rates from the rather weak
condition (\ref{NA1}). Applying a similar reasoning for the
variance of $q_s((X,Y),(X',Y'))$ (which would be the case if one
considered a risk estimate based on independent pairs by splitting
the training data into two halves, see Section \ref{empmin}),
would have led to the condition:
\begin{equation}
\left| \eta (x)-\eta (x^{\prime })\right| \geq c, \,
\end{equation}
for some constant $c$, and $x \neq x^{\prime}$.  This condition is
satisfied only when $\eta (X)$ has a discrete distribution.
\end{remark}

\subsection{Noiseless regression model}

Next we consider the {\em noise-free regression model} in which
$Y=m(X)$ for some (unknown) function $m:\X\to \R$. Here obviously
$L^*=0$ and the Bayes ranking rule is given by the scoring
function $s^*=m$ (or any strictly increasing transformation of
it). Clearly, in this case
\[
  q_r(x,x') = \IND{(m(x)-m(x'))\cdot r(x,x')<0}
\]
and therefore
\[
  \var(h_r(X,Y)) \le \EXP q_r^2(X,X') = L(r)~,
\]
and therefore the condition of Proposition \ref{fast} is satisfied
with $c=1$ and $\alpha=1$. Thus, the risk of the empirical risk
minimizer $r_n$ satisfies, with probability at least $1-\delta$,
\[
L(r_n)  \leq 2 \inf_{r \in \cR} L(r)
+ C \frac{V\log (n/\delta)}{n}
\]
provided $\cR$ has finite {\sc vc} dimension $V$.

\subsection{Regression model with noise}

Now we turn to the {\em general regression model with
heteroscedastic errors} in which $Y=m(X)+ \sigma(X)\epsilon$ for
some (unknown) functions $m:\X\to \R$ and $\sigma:\X\to \R$, where
$\epsilon$ is a standard gaussian random variable, independent of
$X$.

We set
\[
\Delta(X,X')= \frac{m(X)-m(X')}{\sqrt{\sigma^2(X)+\sigma^2(X')}}\,
.
\]

We have again $s^*=m$ (or any strictly increasing transformation
of it) and the optimal risk is
\[L^* = \EXP  \Phi \left(- \left|\Delta(X,X')
\right| \right)
\]
where $\Phi$ is the distribution function of the standard gaussian
random variable. The maximal value of $L^*$ is attained when the
regression function $m(x)$ is constant. Furthermore, we have
\[
L(s)-L^* = \EXP \left(\left| 2 \Phi \left(\Delta(X,X')\right) -1
\right| \cdot \IND{(m(x)-m(x'))\cdot (s(x) - s(x'))<0} \right)~.
\]

\medskip

\noindent {\bf Noise assumption.} There exist constants $c>0$ and
$\alpha\in [0,1]$ such that for all $x\in\X$,
\begin{equation}
\label{NA2} \EXP_{X'}( |\Delta (x, X')|^{-\alpha})\leq c \, .
\end{equation}

\begin{proposition}
Under (\ref{NA2}), we have, for all $s\in\F$
\[
\var(h_s(X,Y)) \leq (2\Phi(c)-1) \,  \Lambda(s)^\alpha \, .
\]
\end{proposition}

\proof
By symmetry, we have
\[
\left| 2 \Phi \left(\Delta(X,X')\right) -1 \right| =  2 \Phi
\left( \left|\Delta(X,X')\right|\right) -1 \, .
\]

Then, using the concavity  of the distribution function $\Phi$ on
$\RR_+$, we have, by Jensen's inequality,
\[
\forall x\in\X, \quad  \EXP_{X'} \Phi( |\Delta (x, X')|^{-\alpha})
\leq \Phi( \EXP_{X'} |\Delta (x,
 X')|^{-\alpha}) \leq  \Phi(c) \, ,
\]
where we have used (\ref{NA2}) together with the fact that $\Phi$
is increasing. Now the result follows following the argument given
in the proof of Proposition \ref{noisecond}.
\qed

The preceding noise condition is fulfilled in many cases, as
illustrated by the example below.

\begin{corollary}
Suppose that $m(X)$ has a bounded density and the conditional
variance $\sigma (x) $ is bounded over $\X$. Then the noise
condition (\ref{NA2}) is satisfied for any $\alpha<1$.
\end{corollary}

\begin{remark}
The argument above still holds if we drop the gaussian noise
assumption. Indeed we only need the random variable $\epsilon$ to
have a symmetric density decreasing over $\RR_+$.
\end{remark}

\section{A moment inequality for $U$-processes}
\label{moment}

In this section we establish a general exponential inequality for
$U$-processes. This result is based on
moment inequalities obtained for empirical processes and
Rademacher chaoses in Bousquet, Boucheron, Lugosi, and Massart
\cite{BoBoLuMa04} and generalizes an inequality
due to Arcones and Gin\'e \cite{ArGi94}. We also refer to the
corresponding results obtained for $U$-statistics by Adamczak
\cite{Ada05}, Gin\'e, Latala, and Zinn \cite{GiLaZi01}, and Houdr\'e
and Reynaud-Bouret \cite{HRB03}.

\begin{theorem}
\label{thm:ineq}
Let $X, X_1, ..., X_n$ be i.i.d. random variables
and let $\F$ be a class of kernels. Consider a
degenerate $U$-process $Z$ of order 2 indexed by $\F$,
\[
Z = \sup_{f\in\F} \left| \sum_{i,j} f(X_i, X_j) \right|
\]
where $\EXP f(X,x)=0$, $\forall x, f$. Assume also
$f(x,x)=0$, $\forall x$ and $\sup_{f\in \F} \|f\|_{\infty} =F$.
 Let $\epsilon_1, ..., \epsilon_n$ be
i.i.d.\ Rademacher random variables and introduce the random variables
\begin{align*}
Z_{\epsilon} & =  \sup_{f\in\F} \left| \sum_{i,j} \epsilon_i
\epsilon_j f(X_i, X_j)  \right|~,\\
U_{\epsilon} &
= \sup_{f\in\F} \sup_{\alpha : \|\alpha\|_2\le 1} \sum_{i,j} \epsilon_i  \alpha_j f(X_i, X_j) ~,\\
M & = \sup_{f\in\F, k=1\ldots n} \left| \sum_{i=1}^n \epsilon_i
f(X_i,X_k) \right|~.
\end{align*}
Then there exists a universal constant $C>0$ such that for
all $n$ and $q\ge 2$,
\[
\left( \EXP Z^q \right)^{1/q}
\le
C \left( \EXP Z_{\epsilon}
+  q^{1/2} \EXP U_{\epsilon}
+ q (\EXP M + F n)
+ q^{3/2}  F n^{1/2}
+
 q^2 F  ~\right) .
\]
Also, there exists a universal constant $C$ such that for all $n$ and
$t>0$,
\[
\PROB\{ Z > C\EXP Z_\epsilon +t\}
 \le \exp\left(-\frac{1}{C}
\min\left( \left(\frac{t}{\EXP U_\epsilon}\right)^2,
                 \frac{t}{\EXP M + Fn},
                 \left(\frac{t}{F\sqrt{n}}\right)^{2/3},
                 \sqrt{\frac{t}{F}} \right) \right)~.
\]
\end{theorem}

\begin{remark}
A generously overestimated value of the constants may be easily deduced
from the proof. We are convinced that these are far from being
the best possible but do not have a good guess of what the best
constants might be.
\end{remark}

\proof
The proof of Theorem \ref{thm:ineq} is based on symmetrization,
decoupling, and concentration inequalities for empirical processes
and Rademacher chaos.

Since the $f$ are degenerate kernels, one may relate the moments of
$Z$ to those of $Z_\epsilon$ by the randomization inequality
%GL check constant! (You had an extra factor of 12 before...
\[
  \EXP Z^q \le 4^{q} \EXP Z_\epsilon^q~,
\]
valid for $q \ge 1$, see Chapter 3 of \cite{deGi99}.
Thus, it suffices to derive moment inequalities for the
symmetrized $U$-process $Z_\epsilon$. We do this by conditioning.
Denote by $\EXP_\epsilon$ the expectation taken with respect to
the variables $\epsilon_i$ (i.e., conditional expectation given
$X_1,\ldots,X_n$). Then we write
$\EXP Z_\epsilon^q = \EXP \EXP_\epsilon Z_\epsilon^q$ and study the
quantity $\EXP_\epsilon Z_\epsilon^q$, with the $X_i$ fixed.
But then $Z_\epsilon$ is a so-called {\sl Rademacher chaos} whose
tail behavior has been studied, see Talagrand \cite{Tal96c},
Ledoux \cite{Led96a}, Boucheron, Bousquet, Lugosi, and Massart
\cite{BoBoLuMa04}.
In particular,
for any $q \ge 2$,
\begin{equation*}
\begin{array}{cll}
\bigl(\EXP_{\epsilon} Z_{\epsilon}^q \bigr)^{1/q} & \le
\EXP_{\epsilon} Z_{\epsilon} + \bigl(\EXP_{\epsilon}
\bigl(Z_{\epsilon}- \EXP_{\epsilon} Z_{\epsilon}\bigr)_+^q
\bigr)^{1/q} & \mbox{(since $Z \ge 0$)}\\
\\
& \le \EXP_{\epsilon} Z_{\epsilon} + 3 \sqrt{q}~\EXP_{\epsilon} U_{\epsilon}
+ 4q B &
\end{array}
\end{equation*}
with $U_{\epsilon}$ defined above and
\[
B  = \sup_{f\in\F} \sup_{\alpha, \alpha' : \|\alpha\|_2,
\|\alpha'\|_2\le 1} \left| \sum_{i,j} \alpha_i \alpha'_j f(X_i,X_j)
\right|
\]
where the second inequality follows by Theorem 14 of \cite{BoBoLuMa04}.
Using the inequality $(a+b+c)^q  \le 3^{q-1}(a^q + b^q +c^q)$
valid for $q \ge 2$, $a,b,c>0$, we have
\[
\EXP_{\epsilon} Z_{\epsilon}^q \le 3^{q-1} \left(
\bigl(\EXP_{\epsilon} Z_{\epsilon}\bigr)^q + 3^q q^{q/2}
\bigl(\EXP_{\epsilon} U_{\epsilon}\bigr)^q +  4^q q^q B^q
\right)~.
\]
It remains to derive suitable upper bounds for the expectation
of the three terms on the right-hand side.

\medskip
\noindent
{\bf First term: $\EXP \bigl(\EXP_{\epsilon} Z_{\epsilon}\bigr)^q$}

In order to handle the moments of
 $\EXP_\epsilon Z_\epsilon$, first we note that by
a decoupling inequality in de la Pe\~na and Gin\'e \cite[page 101]{deGi99},
\[
\EXP_{\epsilon} Z_{\epsilon} \le 8
\EXP_{\epsilon} Z'_{\epsilon}
\]
where
\[
Z'_{\epsilon}= \sup_{f\in\F} \left| \sum_{i,j} \epsilon_i \epsilon'_j
f(X_i, X_j) \right|
\]
Here $\epsilon'_1,\ldots,\epsilon'_n$ are i.i.d.\ Rademacher variables,
independent of the $X_i$ and the $\epsilon_i$.
Nothe that $\EXP_{\epsilon}$ now denotes expectation taken with respect to
both the $\epsilon_i$ and the $\epsilon'_i$.

Thus, we have
\[
\EXP \bigl(\EXP_{\epsilon} Z_{\epsilon}\bigr)^q \le 8^q \EXP
\bigl(\EXP_{\epsilon} Z'_{\epsilon}\bigr)^q
\]
In order to bound the moments of the random variable
$A=\EXP_{\epsilon} Z'_{\epsilon}$, we apply  Corollary 3
of \cite{BoBoLuMa04}.
In order to apply this corollary, define, for $k=1,\ldots,n$,
the random variables
\[
A_k = \EXP_\epsilon
\sup_{f\in\F} \left| \sum_{i,j \neq k} \epsilon_i \epsilon'_j
f(X_i,X_j) \right|
\]
It is easy to see that $A_k \le A$.
%% Moreover, note that, if $f_k^*$ denotes the (random) function that
%% achieves the supremum within the expectation in the definition of
%% $Z_k$, then
%% \[
%%   Z- Z_k \ge \EXP_\epsilon \epsilon_k \sum_{j=1}^n
%%           \epsilon_j' f_k^*(X_k,X_j') = 0
%% \]
%% where the equality follows from the fact that $\epsilon_k$ is
%% independent of the rest of the variables (because of
%% decoupledness!). Similarly, $Z-Z_k' \ge 0$.

On the other hand, defining
\[
    R_k =
\sup_{f\in\F} \left| \sum_{i=1}^n \epsilon_i f(X_i,X_k) \right|~,
\]
we clearly have
\[
   A - A_k \le 2\EXP_\epsilon  R_k~.
\]
Also, denoting by  $f^*$  the (random) function achieving the maximum in
the definition of $Z$, we have
\begin{eqnarray*}
\sum_{k=1}^n (A- A_k)
& \le & \EXP_\epsilon \left( \sum_{k=1}^n \epsilon_k
\sum_{j=1}^n
          \epsilon_j' f^*(X_k,X_j')
+ \sum_{k=1}^n \epsilon_k' \sum_{i=1}^n
          \epsilon_i f^*(X_i,X_k')   \right)   \\*
& = & 2A~,
\end{eqnarray*}
Therefore,
\[
\sum_{k=1}^n (A-A_k)^2 \le 4 A \EXP_\epsilon M
\]
where $ M= \max_k R_k$. Then by Corollary 3 of \cite{BoBoLuMa04},
we obtain
\[
\EXP \bigl(\EXP_{\epsilon} Z'_{\epsilon}\bigr)^q
= \EXP A^q
\le 2^{q-1}
\left( 2^q \bigl(\EXP Z'_{\epsilon}\bigr)^q
  + 5^q q^q \EXP\left(\EXP_\epsilon M\right)^q\right)~.
\]
By un-decoupling (see de la Pe\~na and Gin\'e \cite[page 101]{deGi99}),
we have $\EXP Z'_{\epsilon} \le 4 \EXP Z_{\epsilon}$.

To bound $\EXP\left(\EXP_\epsilon M\right)^q$,
 observe that $\EXP_\epsilon M$ is
a conditional Rademacher average, for which Theorem 13 of
of \cite{BoBoLuMa04} may be applied. According to this,
\[
 \EXP \left(\EXP_\epsilon M\right)^q \le 2^{q-1} \left(
2^q \left( \EXP M \right)^q
+ 5^qq^q F^q
\right)
\]
Collecting terms, we have
\[
\EXP \bigl(\EXP_{\epsilon} Z_{\epsilon}\bigr)^q \le
128^q \left(\EXP Z_{\epsilon}\right)^q  +
320^q q^q \left(\EXP M \right)^q
+  800^q F^q q^{2q}~.
\]

\bigskip

\noindent
{\bf Second term: $\EXP_X \bigl(\EXP_{\epsilon}
U_{\epsilon}\bigr)^q$}

The moments of $\EXP_{\epsilon} U_{\epsilon}$ can be estimated
by the same inequality as the one we used for $\EXP_\epsilon M$
since
$\EXP_{\epsilon} U_{\epsilon}$ is also a conditional Rademacher
average.
Observing that
\[
\sup_{f,i} \sup_{\alpha : \|\alpha\|_2\le 1}
\sum_{j\neq i} \alpha_j f(X_i, X_j)
\le F\sqrt{n}
\]
by the Cauchy-Schwarz inequality, we have, by Theorem 13 from
\cite{BoBoLuMa04},
\[
\EXP \bigl(\EXP_{\epsilon} U_{\epsilon}\bigr)^q
\le 2^{q-1}
 \left( 2^q \bigl(\EXP U_{\epsilon}\bigr)^q + 5^q q^q F^q
n^{q/2} \right)~.
\]

\bigskip

\noindent
{\bf Third term: $\EXP_X B^q$}

Finally, by the Cauchy-Schwarz inequality, we have $B \le nF$ so
\[
\EXP_X B^q \le n^q F^q~.
\]

Now it remains to simply put the pieces together to obtain
\begin{eqnarray*}
\EXP Z^q & \le & 12^q \left( 128^q \bigl(\EXP Z_{\epsilon}\bigr)^q
+ 12^q q^{q/2} \bigl(\EXP U_{\epsilon}\bigr)^q + 320^q q^q
\bigl(\EXP M\bigr)^q + 4^q F^q n^q q^q  \right. \\* & & \left. +
30^q F^q n^{q/2} q^{3q/2} + 800^q F^q q^{2q} \right) ~,
\end{eqnarray*}
proving the announced moment inequality.

In order to derive the exponential inequality, use Markov's
inequality $\PROB\{Z>t\} \le t^{-q} \EXP Z^q$ and choose
\[
q = C \min\left( \left(\frac{t}{\EXP U_\epsilon}\right)^2,
                 \frac{t}{\EXP M},
                 \frac{t}{Fn},
                 \left(\frac{t}{F\sqrt{n}}\right)^{2/3},
                 \sqrt{\frac{t}{F}} \right)
\]
for an appropriate constant $C$.
\qed

\section{Convex risk minimization}\label{convex}

Several successful algorithms for classification, including
various versions of {\em boosting} and {\em support vector
machines} are based on replacing the loss function by a convex
function and minimizing the corresponding empirical convex risk
functionals over a certain class of functions (typically over a
ball in an appropriately chosen Hilbert or Banach space of
functions). This approach has important computational advantages,
as the minimization of the empirical convex functional is often
computationally feasible by gradient descent algorithms. Recently
significant theoretical advance has been made in understanding the
statistical behavior of such methods, see, e.g., Bartlett, Jordan,
and McAuliffe \cite{BaJoMc03}, Blanchard, Lugosi and Vayatis
\cite{BlLuVa03}, Breiman \cite{Bre04}, Jiang \cite{Jia04}, Lugosi and
Vayatis \cite{LuVa04}, Zhang \cite{Zha04}.

The purpose of this section is to extend the principle of convex
risk minimization to the ranking problem studied in this paper.
Our analysis also provides a theoretical framework for the
analysis of some successful ranking algorithms such as the {\sc
RankBoost} algorithm of Freund, Iyer, Schapire, and Singer
\cite{FrIyScSi04}. In what follows we adapt the arguments of Lugosi
and Vayatis \cite{LuVa04} (where a simple binary classification
problem was considered) to the ranking problem.

The basic idea is to consider ranking rules induced by real-valued
functions, that is, ranking rules of the form
\[
  r(x,x') = \left\{ \begin{array}{ll}
                   1 & \mbox{if $f(x,x')>0$}  \\
                  -1 & \mbox{otherwise}
           \end{array} \right.
\]
where $f:\X \times \X \to \R$ is some measurable real-valued
function. With a slight abuse of notation, we will denote by $L(f)
= \PROB\{ \sgn(Z) \cdot f(X,X') < 0\} = L(r)$ the risk of the
ranking rule induced by $f$. (Here $\sgn(x) = 1$ if $x>0$,
$\sgn(x) = -1$ if $x<0$, and $\sgn(x) = 0$ if $x=0$.) Let $\phi:
\R \to [0,\infty)$ be a convex {\em cost function} satisfying
$\phi(0)=1$ and $\phi(x) \ge \IND{x\ge 0}$. Typical choices of
$\phi$ include the exponential cost function $\phi(x)=e^x$, the
"logit" function $\phi(x)=\log_2(1+e^x)$, or the "hinge loss"
$\phi(x)= (1+x)_+$. Define the {\em cost functional} associated to
the cost function $\phi$ by
\[
   A(f) = \EXP \phi(-\sgn(Z) \cdot f(X,X'))~.
\]
Obviously, $L(f) \le A(f)$. We denote by $A^*= \inf_f A(f)$ the
"optimal" value of the cost functional where the infimum is taken
over all measurable functions $f: \X \times \X \to \R$.

The most natural estimate of the cost functional $A(f)$, based on
the training data $D_n$, is the {\em empirical cost functional}
defined by the $U$-statistic
\[
A_n(f) = \frac{1}{n(n-1)}\sum_{i\neq j} \phi(-\sgn(Z_{i,j}) \cdot
f(X_i,X_j))~.
\]
The ranking rules based on {\em convex risk minimization} we
consider in this section minimize, over a set $\F$ of real-valued
functions $f:\X \times \X \to \R$, the empirical cost functional
$A_n$, that is, we choose $f_n = \argmin_{f \in \F} A_n(f)$ and
assign the corresponding ranking rule
\[
  r_n(x,x') = \left\{ \begin{array}{ll}
                   1 & \mbox{if $f_n(x,x')>0$}  \\
                  -1 & \mbox{otherwise.}
           \end{array} \right.
\]
(Here we assume implicitly that the minimum exists. More
precisely, one may define $f_n$ as any function $f\in \F$
satisfying $A_n(f_n) \le \inf_{f \in \F} A_n(f) + 1/n$.)

By minimizing convex risk functionals, one hopes to make the
excess convex risk $A(f_n)-A^*$ small. This
 is meaningful for ranking
if one can relate the excess convex risk to the excess ranking
risk $L(f_n)-L^*$. This may be done quite generally by recalling a
 result of Bartlett, Jordan, and McAuliffe \cite{BaJoMc03}.
To this end, introduce the functions
\[
  H(\rho) = \inf_{\alpha \in \R}
  \left( \rho \phi(-\alpha) + (1-\rho) \phi(\alpha) \right)
\]
and
\[
  H^-(\rho) = \inf_{\alpha: \alpha(2\rho -1)\le 0}
  \left( \rho \phi(-\alpha) + (1-\rho) \phi(\alpha) \right)~.
\]
Defining $\psi$ over $\R$ by
\[
   \psi(x) =  H^-\left(\frac{1+x}{2}\right) - H^-\left(\frac{1-x}{2}\right)~,
\]
Theorem 3 of \cite{BaJoMc03} implies that for all functions $f:\X
\times\X\to \R$,
\[
   L(f) - L^* \le \psi^{-1}\left(A(f) -A^*\right)
\]
where $\psi^{-1}$ denotes the inverse of $\psi$. Bartlett, Jordan,
and McAuliffe show that, whenever $\phi$ is convex, $\lim_{x\to 0}
\psi^{-1}(x) =0$, so convergence of the excess convex risk to zero
implies that the excess ranking risk also converges to zero.
Moreover, in most interesting cases $\psi^{-1}(x)$ may be bounded,
for $x>0$, by a constant multiple of $\sqrt{x}$ (such as in the
case of exponential or logit cost functions) or even by $x$ (e.g.,
if $\phi(x)=(1+x)_+$ is the so-called {\em hinge loss}).

Thus, to analyze the excess ranking risk $L(f)-L^*$ for convex
risk minimization, it suffices to bound the excess convex risk.
This may be done by decomposing it into "estimation" and
"approximation" errors as follows:
\[
   A(f_n) - A^*(f) \le \left( A(f_n) - \inf_{f\in \F} A(f) \right)
      + \left( \inf_{f\in \F} A(f) - A^* \right)~.
\]

Clearly, just like in Section \ref{empmin}, we may (loosely) bound
the excess convex risk over the class $\F$ as
\[
    A(f_n) - \inf_{f\in \F} A(f) \le 2\sup_{f\in \F} |A_n(f)-A(f)|~.
\]
To bound the right-hand side, assume, for simplicity, that the
class $\F$ of functions is uniformly bounded, say
$\sup_{f\in\F,x\in \X} |f(x)|\le B$. Then once again, we may
appeal to Lemma \ref{utoemp} (see the Appendix) and the bounded differences
inequality which imply that for any $\lambda >0$,
\begin{eqnarray*}
\lefteqn{
 \EXP \exp\left(\lambda \sup_{f\in \F} |A_n(f)-A(f)|  \right)
 } \\*
& \le &
 \EXP \exp\left(\lambda \sup_{f\in \F} \left(\frac{1}{\lfloor n/2 \rfloor}
\sum_{i=1}^{\lfloor n/2 \rfloor} \phi\left(-\sgn(Z_{i,\lfloor n/2
\rfloor+i}) \cdot f(X_i,X_{\lfloor n/2 \rfloor+i}) \right) - A(f)
\right)  \right)    \\* & \le &
  \exp\left(\lambda \EXP \sup_{f\in \F} \left(\frac{1}{\lfloor n/2 \rfloor}
\sum_{i=1}^{\lfloor n/2 \rfloor} \phi\left(-\sgn(Z_{i,\lfloor n/2
\rfloor+i}) \cdot f(X_i,X_{\lfloor n/2 \rfloor+i}) \right) - A(f)
\right)  + \frac{\lambda^2 B^2}{2n} \right)~.
\end{eqnarray*}
Now it suffices to derive an upper bound for the expected supremum
appearing in the exponent. This may be done by standard
symmetrization and contraction inequalities. In fact, by mimicking
Koltchinskii and Panchenko \cite{KoPa02} (see also the proof of
Lemma 2 in Lugosi and Vayatis \cite{LuVa04}), we obtain
\begin{eqnarray*}
\EXP \sup_{f\in \F} \left(\frac{1}{\lfloor n/2 \rfloor}
\sum_{i=1}^{\lfloor n/2 \rfloor} \phi\left(-\sgn(Z_{i,\lfloor n/2
\rfloor+i}) \cdot f(X_i,X_{\lfloor n/2 \rfloor+i}) \right) - A(f)
\right) \\*
\le  4B\phi'(B) \EXP \sup_{f\in \F}
\left(\frac{1}{\lfloor n/2 \rfloor} \sum_{i=1}^{\lfloor n/2
\rfloor} \sigma_i  \cdot f(X_i,X_{\lfloor n/2 \rfloor+i}) \right)
\end{eqnarray*}
where $\sigma_1,\ldots,\sigma_{\lfloor n/2 \rfloor}$ i.i.d.\
Rademacher random variables independent of $D_n$, that is,
symmetric sign variables with
$\PROB\{\sigma_i=1\}=\PROB\{\sigma_i=-1\}=1/2$.

We summarize our findings:

\begin{proposition}
\label{rademacherbound} Let $f_n$ be the ranking rule minimizing
the empirical convex risk functional $A_n(f)$ over a class of
functions $f$ uniformly bounded by $-B$ and $B$. Then, with
probability at least $1-\delta$,
\[
A(f_n) - \inf_{f\in \F} A(f) \le 8B\phi'(B) R_n(\F) +
\sqrt{\frac{2B^2\log(1/\delta)}{n}}
\]
where $R_n$ denotes the Rademacher average
\[
R_n(\F) = \EXP \sup_{f\in \F} \left(\frac{1}{\lfloor n/2 \rfloor}
\sum_{i=1}^{\lfloor n/2 \rfloor} \sigma_i  \cdot f(X_i,X_{\lfloor
n/2 \rfloor+i}) \right)~.
\]
\end{proposition}

Many interesting bounds are available for the Rademacher average
of various classes of functions. For example, in analogy of
boosting-type classification problems, one may consider a class
$\F_B$ of functions defined by
\[
 \F_B = \left\{
 f(x,x') = \sum_{j=1}^{N}
w_j g_j(x,x') \, : \,  N\in\NN, \, , \, \sum_{j=1}^{N} |w_j| =B,
\, g_j \in \cR \, \right\}
\]
where $\cR$ is a class of ranking rules as defined in Section
\ref{empmin}. In this case it is easy to see that
\[
   R_n(\F_B) \le B R_n(\cR) \le \mbox{const.}\frac{BV}{\sqrt{n}}
\]
where $V$ is the {\sc vc} dimension of the "base" class $\cR$.

Summarizing, we have shown that a ranking rule based on the
empirical minimization $A_n(f)$ over a class of ranking functions
$\F_B$ of the form defined above, the excess ranking risk
satisfies, with probability at least $1-\delta$,
\[
  L(f_n) - L^* \le \psi^{-1} \left(
8B\phi'(B) c \frac{BV}{\sqrt{n}}
 + \sqrt{\frac{2B^2\log(1/\delta)}{n}}
+  \left( \inf_{f\in \F_B} A(f) - A^* \right) \right)~.
\]
This inequality may be used to derive the {\em universal
consistency} of such ranking rules. For example, the following
corollary is immediate.

\begin{corollary}
Let $\cR$ be a class of ranking rules of finite {\sc vc} dimension
$V$ such that the associated class of functions $\F_B$ is rich in
the sense that
\[
\lim_{B\to \infty} \inf_{f\in \F_B} A(f) = A^*
\]
for all distributions of $(X,Y)$. Then if $f_n$ is defined as the
empirical minimizer of $A_n(f)$ over $\F_{B_n}$ where the sequence
$B_n$ satisfies $B_n \to \infty$ and
% NV: rate corrected here.
$B_n^2 \phi'(B_n)/\sqrt{n} \to 0$, then
\[
\lim_{n\to \infty} L(f_n) = L^* \quad \mbox{almost surely}.
\]
\end{corollary}

Classes $\cR$ satisfying the conditions of the corollary exist, we
refer the reader to Lugosi and Vayatis \cite{LuVa04} for several
examples.

Proposition \ref{rademacherbound} can also be used for
establishing performance bounds for kernel methods such as support
vector machines. A prototypical kernel-based ranking method may be
defined as follows. To lighten notation, we write ${\cal W} =
\X\times\X$.

Let $k:{\cal W}\times{\cal W}\to \R$ be a symmetric positive
definite function, that is,
\[
\sum_{i,j=1}^n \alpha_i\alpha_j k(w_i,w_j)\ge 0\,,
\]
for all choices of $n$, $\alpha_1,\ldots,\alpha_n\in\R$ and
$w_1,\ldots,w_n\in {\cal W}$.

A kernel-type ranking algorithm may be defined as one that
performs minimization of the empirical convex risk $A_n(f)$
(typically based on the hinge loss $\phi(x)=(1+x)_+$) over the
class $\F_{B}$ of functions defined by
 a ball of the associated reproducing kernel
Hilbert space of the form (where $w=(x,x')$)
\[
  \F_{B} = \left\{ f(w) = \sum_{j=1}^N c_j k(w_j,w): \ N \in \mathbb{N}, \
   \sum_{i,j=1}^N c_ic_jk(w_i,w_j) \le B^2, w_1,\ldots,w_N \in {\cal W}
   \right\}\,.
\]
In this case we have
\[
R_n(\F_B) \le \frac{2B}{n} \EXP \sqrt{\sum_{i=1}^{\lfloor n/2
\rfloor} k((X_i,X_{\lfloor n/2 \rfloor+i}),(X_i,X_{\lfloor n/2
\rfloor+i}))}~,
\]
see, for example, Boucheron, Bousquet, and Lugosi \cite{BoBoLu05}.
Once again, universal consistency of such kernel-based ranking
rules may be derived in a straightforward way if the approximation
error $\inf_{f\in \F_B} A(f) - A^*$ can be guaranteed to go to
zero as $B\to \infty$. For the approximation properties of such
kernel classes we refer the reader to Cucker and Smale
\cite{CuSm02}, Scovel and Steinwart
 \cite{ScSt03}, Smale and Zhou \cite{SmZh03}, Steinwart \cite{Ste01}, etc.

\section*{Appendix 1: Basic facts about $U$-statistics}

Here we recall some basic facts about $U$-statistics.
Consider the i.i.d. random variables $X, X_1, ..., X_n$  and
denote by
\[
U_n = \frac{1}{n(n-1)} \sum_{i \neq j} q(X_i, X_j)
\]
a $U$-statistic of order $2$ where the kernel $q$ is a symmetric
real-valued function.

$U$-statistics have been studied in depth and their behavior is
well understood. One of the classical inequalities concerning
$U$-statistics is due to Hoeffding \cite{Hoe63} which implies
that, for all $t>0$,
\[
   \PROB\{ |U_n - \EXP U_n| > t\} \le 2e^{-2\lfloor (n/2) \rfloor t^2}
                          \le 2e^{-(n-1)t^2}~.
\]

Hoeffding also shows that, if $\sigma^2 = \var(q(X_1,X_2))$, then
\begin{equation}
\label{bernstein}
   \PROB\{ |U_n - \EXP U_n| > t\}
   \le 2\exp\left(-\frac{\lfloor (n/2) \rfloor t^2}{2\sigma^2 + 2t/3}\right)~.
\end{equation}

It is important noticing here that the latter inequality may be
improved by replacing $\sigma^2$ by a smaller term. This is based
on the so-called Hoeffding's decomposition as described below.

\medskip

 The
$U$-statistic $U_n$ is said {\em degenerate} if its kernel $q$
satisfies
\[
\forall x, \quad \EXP \left( q(x,X) \right) =0~.
\]

There are two basic representations of $U$-statistics which we
recall next (see Serfling \cite{Ser80} for more details).

\medskip

\noindent {\bf Average of 'sums-of-i.i.d.' blocks}

\medskip

This representation is the key for obtaining 'first-order' results
for non-degenerate $U$-statistics. The $U$-statistic $U_n$ can be
expressed as
\[
U_n = \frac{1}{n!} \sum_{\pi} \frac{1}{\lfloor n/2
\rfloor} \sum_{i=1}^{\lfloor n/2 \rfloor}
q\bigl(X_{\pi(i)},X_{\pi(\lfloor n/2 \rfloor+i)}\bigr)
\]
where the sum is taken over all permutations $\pi$  of
$\{1,\ldots,n\}$. The idea underlying this representation is to
reduce the analysis to the case of sums of i.i.d. random
variables. The next simple lemma is based on this representation.

\begin{lemma}
\label{utoemp} Let $q_\tau: \X\times\X \to \R$ be real-valued
functions indexed by $\tau \in T$ where $T$ is some set. If
$X_1,\ldots,X_n$ are i.i.d.\ then for any convex nondecreasing
function $\psi$,
\begin{eqnarray*}
\lefteqn{
  \EXP \psi\left( \sup_{\tau \in T} \frac{1}{n(n-1)}\sum_{i\neq j}
      q_{\tau}(X_i,X_j) \right)  } \\*
& \le &
  \EXP \psi\left( \sup_{\tau \in T} \frac{1}{\lfloor n/2 \rfloor}
\sum_{i=1}^{\lfloor n/2 \rfloor} q_\tau(X_i,X_{\lfloor n/2
\rfloor+i})
      \right)~,
\end{eqnarray*}
assuming the suprema are measurable and the expected values exist.
\end{lemma}

\proof The proof uses the same trick Hoeffding's above-mentioned
inequalities are based on. Observe that
\begin{eqnarray*}
\lefteqn{
  \EXP \psi\left( \sup_{\tau \in T} \frac{1}{n(n-1)}\sum_{i\neq j}
      q_{\tau}(X_i,X_j) \right)    } \\*
& = &
  \EXP \psi\left( \sup_{\tau \in T}
\frac{1}{n!} \sum_{\pi} \frac{1}{\lfloor n/2 \rfloor}
\sum_{i=1}^{\lfloor n/2 \rfloor} q_\tau(X_{\pi(i)},X_{\pi(\lfloor
n/2 \rfloor+i)})
   \right) \\
& \le &
  \EXP \psi\left(
\frac{1}{n!} \sum_{\pi} \sup_{\tau \in T} \frac{1}{\lfloor n/2
\rfloor} \sum_{i=1}^{\lfloor n/2 \rfloor}
q_\tau(X_{\pi(i)},X_{\pi(\lfloor n/2 \rfloor+i)})
   \right) \\*
& & \quad \mbox{(since $\psi$ is non-decreasing)}\\
& \le & \frac{1}{n!} \sum_{\pi}  \EXP \psi\left(
 \sup_{\tau \in T} \frac{1}{\lfloor n/2 \rfloor}
\sum_{i=1}^{\lfloor n/2 \rfloor} q_\tau(X_{\pi(i)},X_{\pi(\lfloor
n/2 \rfloor+i)})
   \right) \\*
& & \quad \mbox{(by Jensen's inequality)}\\
& = &
  \EXP \psi\left( \sup_{\tau \in T} \frac{1}{\lfloor n/2 \rfloor}
\sum_{i=1}^{\lfloor n/2 \rfloor} q_\tau(X_i,X_{\lfloor n/2
\rfloor+i})
      \right)
\end{eqnarray*}
as desired. \qed

\bigskip

\noindent {\bf Hoeffding's decomposition}
%NV -historical note:
%    Hoeffding's decomposition was apparently given in a confidential TR
%    by Hoeffding in 1961, then it was independently obtained by Berk in 1966.
%    A historical account on this story is given by  Oosterhoff and van Zwet
%    in the book edited by Fisher and Sen ``The collected works of Wassily
%    Hoeffding'' (Springer).
%    Otherwise, a classical reference for U-statistics is the book by Serfling
%    ``Approximation Theorems of Mathematical Statistics''.

\medskip

Another way to interpret a $U$-statistics is as an orthogonal
expansion known as Hoeffding's decomposition.

Assuming that $q(X_1, X_2)$ is square integrable, $U_n-\EXP U_n$
may be decomposed as a sum $T_n$ of i.i.d. random variables plus a {\em
degenerate} $U$-statistic $W_n$. In order to write this
decomposition, consider the following function of one variable
\[
h(X_i) = \EXP (q(X_i, X) \mid X_i) - \EXP U_n  \, ,
\]
and the function of two variables
\[
\wh{h}(X_i, X_j)= q(X_i, X_j) - \EXP U_n- h(X_i) - h(X_j).
\]
Then we have the orthogonal expansion
\[
U_n = \EXP U_n + 2 T_n + W_n~,
\]
where
\begin{align*}
T_n & = \frac{1}{n} \sum_{i=1}^n h(X_i), \\
W_n &
    = \frac{1}{n(n-1)} \sum_{i \neq j} \wh{h}(X_i, X_j) \, .
\end{align*}
$W_n$ is a degenerate $U$-statistic because its kernel $\wh{h}$
satisfies
\[
   \EXP \left( \wh{h}(X_i,X) \mid X_i\right) =0~.
\]
Clearly, the variance of $T_n$ is
\[
\Var(T_n) = \frac{\Var(\EXP (q(X_1, X) \mid X_1))}{n} \, .
\]
Note that $\Var(\EXP (q(X_1, X) \mid X_1))$ is less than
$\Var(q(X_1, X))$ (unless $q$ is already degenerate). Furthermore,
the variance of the degenerate $U$-statistic $W_n$ is of the order
$1/n^2$. $T_n$ is thus the leading term in this orthogonal
decomposition. Indeed, the limit distribution of  $\sqrt{n}(U_n
-\EXP U_n)$ is the normal distribution $\N(0, 4\Var(\EXP (q(X_1,X)
\mid X_1))$ (see \cite{Hoe48}). This suggests that inequality
(\ref{bernstein}) may be quite loose.

Indeed, exploiting further Hoeffding's decomposition (combined
with arguments related to decoupling, randomization and
hypercontractivity of Ra\-de\-ma\-cher chaos) de la Pe\~na and
Gin\'e \cite{deGi99} established a Bernstein's type inequality
of the form (\ref{bernstein}) but with $\sigma^2$ replaced by the
variance of the conditional expectation (see Theorem 4.1.13 in
\cite{deGi99}).

Specialized to our setting with $q(X_i, X_j) =  \IND{Z_{i,j} \cdot
r(X_i,X_j) <0}$ the inequality of de la Pe\~na and Gin\'e states
that
\[
   \PROB\{ |L_n(r) - L(r)| > t\}
   \le 4\exp\left(-\frac{n t^2}{8s^2 + c t}\right)~,
\]
where $s^2 = \var(\PROB\{Z \cdot r(X,X') <0|X\})$ is the variance
of the conditional expectation and $c$ is some constant.

\section*{Appendix 2: Connection with the {\sc roc} curve and the {\sc auc} criterion}

In the bipartite ranking problem, the {\sc roc} curve ({\sc roc} standing for
\textit{Receiving Operator Characteristic}, see
\cite{GrSw66}) and the {\sc auc} criterion are popular measures for
evaluating the performance of scoring functions in
applications.

\medskip

Let $s:\mathcal{X}\rightarrow \mathbb{R}$ be a scoring function.
The {\sc roc} curve is defined by plotting the \textit{true
positive rate}
\[
\mbox{{\sc tpr}}_s(x) = \PROB\left(s(X\right) \geq x\mid Y=1)
\]
against the \textit{false positive rate}
\[
\mbox{{\sc fpr}}_s(x)=\PROB\left(s(X\right) \geq x\mid Y=-1) \, .
\]

By a straightforward change of parameter, the {\sc roc} curve may be
expressed as the graph of the power of the test defined by $s(X)$
as a function of its level $\alpha$:
\[
\beta_{s}(\alpha) = \mbox{{\sc tpr}}_s(q_{s,\alpha})
\]
where $q_{s,\alpha }
=\inf \{x\in (0,1) \, : \, \mbox{{\sc fpr}}_s(x) \leq \alpha \}$.

\medskip

Observe that if $s(X)$ and $Y$ are independent (i.e.,
when $\mbox{{\sc tpr}}_s=\mbox{{\sc fpr}}_s$), the {\sc roc} curve is
simply the diagonal segment $\beta_s(\alpha)=\alpha$. This
measure of accuracy induces a partial order on the set of all
scoring functions: for any $s_{1},$ $s_{2}$,  we say that
$s_{1}$ is more accurate than $s_{2}$ if and only if its {\sc roc} curve
is above the one of $s_{2}$ for every level $\alpha$, that is,
 if and only if $\beta _{s_{2}}(\alpha )\le \beta
_{s_{1}}(\alpha )$ for all $\alpha \in (0,1)$.

\begin{proposition}\label{optimalroc}
The regression function $\eta$ induces an optimal ordering on
$\mathcal{X}$ in the sense that its {\sc roc} curve is not below
any other scoring function $s$:
\[
\forall\alpha \in [0,1] ,  \quad  \beta _{\eta }(\alpha ) \ge
\beta _{s}(\alpha ) .
\]
\end{proposition}

\proof
The result follows from the Neyman-Pearson lemma
applied to the test of the null assumption "$Y=-1$" against the
alternative "$Y=1$" based on the observation $X$: the test based
on the likelihood ratio $\eta(X)/(1-\eta(X))$ is uniformly more
powerful than any other test based on $X$.
\qed

\begin{remark}
Note that the {\sc roc} curve does not characterize the scoring
function. For any $s$ and any strictly increasing
function $h \, : \, \mathbb{R}\rightarrow \mathbb{R}$, $s$ and
$h\circ s$ clearly yield the same ordering on $\mathcal{X}:$
$\beta _{s}=\beta _{h\circ s}.$
\end{remark}

%\begin{remark}
%Another criterion for evaluating the accuracy of a scoring function that may
%be found in the literature devoted to information retrieval is the
%\textit{precision-recall curve} $x\in \mathbb{R}\mapsto (\PROB(Y=1\mid
%s(X)>x), \PROB\left( s(X\right) >x\mid Y=1))$.
%\end{remark}

Instead of optimizing the {\sc roc} curve over a class of scoring
functions which is a difficult task, a simple idea is to
search for $s$ that maximizes the Area Under the {\sc roc} Curve (known
as the {\sc auc} criterion) :
\[
\mbox{{\sc auc}}(s)=\int_{0}^{1}\beta _{s}(\alpha ) \, d\alpha \, .
\]

This theoretical quantity may be easily interpreted in a
probabilistic fashion as shown by the following proposition.

\begin{proposition}
For any scoring function $s$,
\[
\mbox{{\sc auc}}(s)= \PROB\left( s(X) \geq s(X^{\prime}) \mid Y=1, \,
Y^{\prime}=-1\right) \, ,
\]
where $(X, Y)$ and $(X^{\prime},Y^{\prime })$ are independent
pairs drawn from the binary classification model.
\end{proposition}

\proof Let $U$ be a  uniformly distributed random variable over
$(0,1)$, independent of $(X,Y)$. Denote by $F_s$ the distribution
function of $s(X)$ given $Y=-1$. Then
\begin{align*}
\mbox{{\sc auc}}(s) & = \int_{0}^{1}\PROB\left( s(X) \ge  q_{s,\alpha}
 \mid Y=1\right) \, d\alpha  \notag \\
& = \EXP(\PROB(s(X) \geq F_{s}^{-1}(U) \mid
Y=1))  \notag \\
& = \PROB\left( s(X) \geq s(X^{\prime}) \mid Y=1, \,
Y^{\prime}=-1\right) \, .
\end{align*}
\qed

\bigskip

\noindent{\bf Acknowledgements.} We thank  Gilles Blanchard and
G\'erard Biau for their valuable comments on a previous version of
this manuscript.

%\bibliography{book1,book2,sajat,cikkek}

\begin{thebibliography}{10}

\bibitem{Ada05}
R.~Adamczak.
\newblock Moment inequalities for $U$-statistics.
\newblock Technical report, Institute of Mathematics of the Polish Academy of
  Sciences, 2005.

\bibitem{AgGrHeHaRo05}
S.~Agarwal, T.~Graepel, R.~Herbrich, S.~Har-Peled, and D.~Roth.
\newblock Generalization bounds for the area under the {ROC} curve.
\newblock {\em Journal of Machine Learning Research}, 6:393--425, 2005.

\bibitem{ArGi93}
M.~A. Arcones and E.~Gin{\'e}.
\newblock Limit theorems for {$U$}-processes.
\newblock {\em The Annals of Probability}, 21:1494--1542, 1993.

\bibitem{ArGi94}
M.~A. Arcones and E.~Gin\'e.
\newblock $U$-processes indexed by {V}apnik-{C}ervonenkis classes of functions
  with applications to asymptotics and bootstrap of $u$-statistics with
  estimated parameters.
\newblock {\em Stochastic Processes and their Applications}, 52:17--38, 1994.

\bibitem{BaJoMc03}
P.L. Bartlett, M.~I. Jordan, and J.~D. McAuliffe.
\newblock Convexity, classification, and risk bounds.
\newblock {\em Journal of the American Statistical Association}, 2005.

\bibitem{BaMe06}
P.L. Bartlett and S.~Mendelson.
\newblock Empirical minimization.
\newblock {\em Probability Theory and Related Fields}, 135, 2006.

\bibitem{BlLuVa03}
G.~Blanchard, G.~Lugosi, and N.~Vayatis.
\newblock On the rates of convergence of regularized boosting classifiers.
\newblock {\em Journal of Machine Learning Research}, 4:861--894, 2003.

\bibitem{BoBoLu05}
S.~Boucheron, O.~Bousquet, and G.~Lugosi.
\newblock Theory of classification: a survey of some recent advances.
\newblock {\em ESAIM. Probability and Statistics}, 9:323--375, 2005.

\bibitem{BoBoLuMa04}
S.~Boucheron, O.~Bousquet, G.~Lugosi, and P.~Massart.
\newblock Moment inequalities for functions of independent random variables.
\newblock {\em The Annals Probability}, 33:514--560, 2005.

\bibitem{Bre04}
L.~Breiman.
\newblock  Population theory for boosting ensembles.
\newblock {\em Annals of Statistics}, 32:1--11, 2004.

\bibitem{CuSm02}
F.~Cucker and S.~ Smale.
\newblock On the mathematical foundations of learning.
\newblock  {\em Bulletin of the American Mathematical Society}, 39:1--49, 2002.

\bibitem{deGi99}
V.H. de~la Pe\~na and E.~Gin\'e.
\newblock {\em Decoupling: from Dependence to Independence}.
\newblock Springer, New York, 1999.

\bibitem{DeGyLu95}
L.~Devroye, L.~Gy\"orfi, and G.~Lugosi.
\newblock {\em A Probabilistic Theory of Pattern Recognition}.
\newblock Springer-Verlag, New York, 1996.

\bibitem{FrIyScSi04}
Y.~Freund, R.~Iyer, R.E. Schapire, and Y.~Singer.
\newblock An efficient boosting algorithm for combining preferences.
\newblock {\em Journal of Machine Learning Research}, 4(6):933--969, 2004.

\bibitem{FHT00}  J.~Friedman, T.~Hastie, and R.~Tibshirani.
\newblock Additive logistic regression: a statistical view of boosting (with discussion).
\newblock  {\em Annals of Statistics}, 28:307--337, 2000.

\bibitem{GiLaZi01}
E.~Gin\'e, R.~Lata{\l}a, and J.~Zinn.
\newblock Exponential and moment inequalities for {U}-statistics.
\newblock In {\em High Dimensional Probability II---Progress in Probability},
  pages 13--38. Birkhauser, 2000.

\bibitem{GiZi84}
E.~Gin\'e and J.~Zinn.
\newblock Some limit theorems for empirical processes.
\newblock {\em Annals of Probability}, 12:929--989, 1984.

\bibitem{GrSw66}  D.M.~Green and J.A.~Swets (1966).
\newblock\textit{Signal detection theory and psychophysics.}
\newblock Wiley, NY.

\bibitem{Hau95}
D.~Haussler.
\newblock Sphere packing numbers for subsets of the boolean $n$-cube with
  bounded {V}apnik-{C}hervonenkis dimension.
\newblock {\em Journal of Combinatorial Theory, Series A}, 69:217--232, 1995.

\bibitem{Hoe48}
W.~Hoeffding.
\newblock A class of statistics with asymptotically normal distributions.
\newblock {\em Annals of Mathematical Statistics}, 10:293--325, 1948.

\bibitem{Hoe63}
W.~Hoeffding.
\newblock Probability inequalities for sums of bounded random variables.
\newblock {\em Journal of the American Statistical Association}, 58:13--30,
  1963.

\bibitem{HRB03} C.~Houdr\'e and P.~Reynaud-Bouret.
\newblock Exponential Inequalities, with constants, for
$U$-statistics of order two.
\newblock Stochastic Inequalities and Applications - Progress in Probability,
Birkhauser, 2003.

\bibitem{Jia04}
W.~Jiang.
\newblock Process consistency for Adaboost (with discussion).
\newblock {\em Annals of Statistics}, 32:13--29, 2004.

\bibitem{Kol06}
V.~Koltchinskii.
\newblock Local {R}ademacher complexities and oracle inequalities in risk
  minimization.
\newblock {\em Annals of Statistics}, 36:00--00, 2006.

\bibitem{KoPa02}
V.~Koltchinskii and D.~Panchenko.
\newblock Empirical margin distribution and bounding the generalization error of combined classifiers.
\newblock {\em Annals of Statistics}, 30:1--50, 2002.

\bibitem{Led96a}
M.~Ledoux.
\newblock On {T}alagrand's deviation inequalities for product measures.
\newblock {\em ESAIM: Probability and Statistics}, 1:63--87, 1997.
\newblock {\tt http://www.emath.fr/ps/}.

\bibitem{Lug02}
G.~Lugosi.
\newblock Pattern classification and learning theory.
\newblock In L.~Gy\"orfi, editor, {\em Principles of Nonparametric Learning},
  pages 5--62. Springer, Wien, 2002.

\bibitem{LuVa04}
G.~Lugosi and N.~Vayatis.
\newblock On the {B}ayes-risk consistency of regularized boosting methods (with discussion).
\newblock {\em Annals of Statistics}, 32:30--55, 2004.

\bibitem{Mas06}
P.~Massart.
\newblock {\em Concentration inequalities and model selection}.
\newblock Ecole d'\'et\'e de Probabilit\'es de Saint-Flour 2003. Lecture Notes
  in Mathematics. Springer, 2006.

\bibitem{MaNe06}
P.~Massart and E.~N\'ed\'elec.
\newblock Risk bounds for statistical learning.
\newblock {\em Annals of Statistics}, 34, 2006.

\bibitem{McD89}
C.~McDiarmid.
\newblock On the method of bounded differences.
\newblock In {\em Surveys in Combinatorics 1989}, pages 148--188. Cambridge
  University Press, Cambridge, 1989.

\bibitem{ScSt03}
S.~Scovel and I.~Steinwart.
\newblock Fast Rates for Support Vector Machines.
\newblock Technical Report LA-UR-03-9117, Los Alamos National Laboratory, 2003.


\bibitem{Ser80} R.J.~Serfling.
\newblock Approximation theorems of mathematical statistics.
\newblock John Wiley \& Sons, 1980.

\bibitem{SmZh03}  S.~Smale and D.-X.~Zhou.
\newblock Estimating the approximation error in learning theory.
\newblock Analysis and Applications, 1, pp. 17-41. Support Vector Machine Soft Margin Classifiers, 2003.

\bibitem{Ste01}
I.~Steinwart.
\newblock On the influence of the kernel on the consistency of support vector machines.
\newblock {\em Journal of Machine Learning Research}, 2:67--93, 2001.

\bibitem{Stu91}
W.~Stute.
\newblock Conditional $U$-statistics.
\newblock {\em Annals of Probability}, 19:812--825, 1991.

\bibitem{Stu94}
W.~Stute.
\newblock Universally consistent conditional {$U$}-statistics.
\newblock {\em The Annals of Statistics}, 22:460--473, 1994.

\bibitem{Tal96c}
M.~Talagrand.
\newblock New concentration inequalities in product spaces.
\newblock {\em Inventiones Mathematicae}, 126:505--563, 1996.

\bibitem{Tsy02}
A.~B. Tsybakov.
\newblock Optimal aggregation of classifiers in statistical learning.
\newblock {\em Annals of Statistics}, 32:135--166, 2004.

\bibitem{VaCh74a}
V.N. Vapnik and A.Ya. Chervonenkis.
\newblock {\em Theory of Pattern Recognition}.
\newblock Nauka, Moscow, 1974.
\newblock (in Russian); German translation: {\em Theorie der Zeichenerkennung},
  Akademie Verlag, Berlin, 1979.

\bibitem{Zha04}
T.~Zhang.
\newblock Statistical behavior and consistency of classification methods based on convex risk minimization (with discussion).
\newblock Annals of Statistics, 32:56--85, 2004.


\end{thebibliography}

\end{document}